%% file: imagedata_submitted.tex
\documentclass{amsart}

\usepackage{graphicx}
\usepackage{color}
\usepackage{color}
\definecolor{newblue}{RGB}{94,89,144}
\definecolor{newblue2}{cmyk}{1,0.6,0,0.06}
\definecolor{orange}{RGB}{200,90,0}
\usepackage[pdfborder={0 0 0},
  colorlinks=true,
  citecolor=newblue,
  linkcolor=newblue2,
  urlcolor=newblue]{hyperref}
\usepackage{subcaption}

\def\mbR{\mathbb{R}}
\newcommand{\mtd}[1]{\partial_t^{(#1)}}

\title[RDS on evolving surfaces from images]{Solving reaction-diffusion equations on evolving surfaces defined by biological image data}
\author[T Bretschneider]{Till Bretschneider}
\author[CJ Du]{Cheng-Jin Du}
\author[CM Elliott]{Charles M. Elliott}
\author[T Ranner]{Thomas Ranner}
\author[B Stinner]{Bj\"{o}rn Stinner}
\address[T Bretschneider, CJ Du]{Warwick Systems Biology Centre, Senate House, University of Warwick, Coventry CV4 7AL, United Kingdom}
\address[CM Elliott,B Stinner]{Mathematics Institute, Zeeman Building, University of Warwick, Coventry CV4 7AL, United Kingdom}
\address[T Ranner]{School of Computing, EC Stoner Building, University of Leeds, LS2 9JT Leeds, United Kingdom}
\email{T.Bretschneider@warwick.ac.uk, C.Du@warwick.ac.uk, C.M.Elliott@warwick.ac.uk, T.Ranner@leeds.ac.uk, Bjorn.Stinner@warwick.ac.uk}

\begin{document}

\begin{abstract} 
We present a computational approach for solving reaction-diffusion equations on evolving surfaces which have been 
obtained from cell image data. It is based on  finite element spaces  defined on  surface triangulations extracted  from time series of 3D images. 
A model for the transport of material between the subsequent surfaces is required where we postulate a velocity in normal direction. 
We apply the technique to image data obtained from a spreading neutrophil cell. By simulating FRAP experiments we investigate the impact of the evolving geometry on the recovery. We find that for idealised FRAP conditions, changes in membrane geometry, easily account for differences of $\times 10$ in recovery half-times, which shows that experimentalists must take great care when interpreting membrane photobleaching results. We also numerically solve an activator -- depleted substrate system 
and report on the effect of the membrane movement on the pattern evolution.
\end{abstract}


\maketitle


\input{source/introduction_submitted.tex}


\input{source/methods_submitted.tex}


\input{source/application_submitted.tex}


\input{source/conclusion_submitted.tex}


\bibliographystyle{plain}
\bibliography{./bib/library}

\end{document}

%% file: source/introduction_submitted.tex
\section{Introduction}

Progress in cell microscopy has led to the availability in large quantities of time dependent high resolution images of cells. These may be processed to yield
details of evolving cell shapes. Furthermore via fluorescent tagging constituents of the cell membrane, lipids or proteins, can 
be visualised. This leads to the possibility of assimilating this data within mathematical models for cell signalling, for example in the context of cell migration,
and also to develop, investigate and test new mathematical models. 

In this article we take the cell to be an evolving three-dimensional domain and we are concerned with partial differential equation (PDE) models for processes 
which take place on the cell membrane. The general goal is to study the effect of the surface evolution on the solution behaviour. We present a numerical approach 
to solve the PDEs on surfaces that are obtained from postprocessing the image data. This image-based modelling is in the spirit 
of \cite{Sba13} but now taking into account time dependent domains 
obtained from  observed data.

The observations comprise a time series of images. Each item of the time series has a selection of two dimensional cross-sections (slices) of the three dimensional object. 
Each slice is segmented in order to identify the cell and its boundary. The upshot is a time series of 
triangulated genus zero surfaces defined by the three dimensional  coordinates of triangle vertices.  Spherical parameterisation using an equiareal mapping to the unit sphere \cite{BreGerKue1995} together 
with  resampling of the mesh ensure that a mesh of constant topology and good quality is maintained. One can then follow the discrete path of a vertex by the 
knowledge of the position at discrete times yielding a discrete \emph{vertex velocity} obtained by time differencing the location of vertices at successive time steps. 


It is natural to use finite element technology in order to solve the PDE models on the triangulations which describe the evolution.
Our computational approach is based on an \emph{arbitrary Lagrangian--Eulerian} (ALE) \emph{evolving surface finite element method} (ESFEM) \cite{DziEll07,EllSty12,EllVen15}. 
The underlying continuum problem formulation requires a given evolving surface including a \emph{material velocity} field related to the trajectories of material surface points. 
Our notion is that it describes how the molecules forming the cell boundary would move without any reactions and lateral transport such as diffusion. 
The contribution of the material velocity in the normal direction describes the geometric evolution of the surface. In general, it may also have a tangential contribution 
which is associated with lateral transport of mass.

The spatial discretisation of the continuum problem involves approximating the evolving surface by an evolving triangulated surface, and approximating the material velocity by a 
\emph{discrete material velocity} which is defined on the evolving triangulated surface.
In the present context the evolving triangulated surface is known only at discrete times given by the image data acquisition. No information as to 
where material points move is provided. In order to define the discrete material velocity we postulate a material velocity in the normal direction. 
The discrete material velocity is defined on each triangle by taking the normal component of the vertex velocity. This means using the tangential component of the 
vertex velocity as the arbitrary tangential velocity in the ALE approach. 

We briefly mention that there are other approaches to solving PDEs on evolving surfaces using fixed grids which include space time finite elements, \cite{OlsReu14}, implicit surface level set methods, \cite{DziEll08b}, 
a Lagrangian level set and particle approach \cite{BerSbaKou10}, phase field methods, \cite{EllStiSty11} and closest point methods, \cite{PetRut16}.

\subsection{Data accessibility}

The original image data as well as segmented individual slices can be obtained from:
\begin{quote}
  \url{http://www2.warwick.ac.uk/fac/sci/systemsbiology/staff/bretschneider/data/}
\end{quote}
The triangulation data and the code used for the simulations are provided as part of the supplementary material and can be obtained from:
\begin{quote}
 \url{https://github.com/tranner/dune-imagedata}
\end{quote}

\subsection{Competing interests}

We have no competing interests.

\subsection{Author's contributions}

T Bretschneider made the images data accessible and discussed the biological applications. 
CJ Du wrote and applied the post-processing software for the image data.
CM Elliott contributed background information,  was involved in the design of the numerical experiments and drafted the manuscript. 
T Ranner was involved in the design of the numerical experiments and performed and post-processed the simulations. 
B Stinner drafted the manuscript, described the method for the surface partial differential equation and was involved in the design of the numerical experiments.

\subsection{Acknowledgements and Funding}

The authors would like to thank the Isaac Newton Institute for Mathematical Sciences for its hospitality during the 
programme \emph{Coupling Geometric PDEs with Physics for Cell Morphology, Motility and Pattern Formation} supported by EPSRC Grant Number EP/K032208/1. 
This research was partially funded by HEIF5 via the Warwick Impact Fund. 
The research of TR was funded by the Engineering and Physical Sciences Research Council (EPSRC EP/L504993/1).
The research of CME was partially supported by the Royal Society via a Wolfson Research Merit Award.
TB was supported by BBSRC grants BBI0082091 and BB/M01150X/1.

%% file: source/methods_submitted.tex
\section{Models and methods}

\subsection{Image data acquisition and processing}
\label{sec:data_ac}

We briefly sketch the method and refer to \cite{DuHawSte13} for more details. First a downsampled version of each image is segmented
using a graph-based approach 
with automatic seed detection. Secondly, since we know the cell boundary is a 
genus zero surface  any topological artefact generated via the first process must be fixed. For example,  small protrusions which cross at a distance away from 
the body of the  cell may result in  holes in the segmented image. A level set approach is used to smooth any 
protrusions until we have a simply connected surface. From the level set representation a spherical parameterisation is obtained based on 
the method of \cite{BreGerKue1995}.
This is then used to push forward a fixed nice triangulation of the sphere on to the actual surface.

\subsection{Continuum model}

In the continuum model underlying the computational approach, the cell boundary is denoted by $\{\Gamma(t)\}_{t \in [0,T]}$ where 
$\Gamma(t)$ is a topological sphere embedded in $\mathbb{R}^3$ at times $t \in [0,T]$ with a unit normal vector field, pointing out of the 
cell, $\nu(\cdot, t): \Gamma(t) \to \mbR^3$, $t \in [0,T]$. 
We assume that a material particle $p$ located at $x_p(t) \in \Gamma(t)$ at a time $t \in [0,T]$ has a velocity $\dot{x}_p(t)$ with both a
 normal component (which determines the evolution of the shape) and tangential components (which are related to transport of material along the surface). 
 A velocity field $v(\cdot, t) : \Gamma(t) \to \mbR^3$, $t \in (0,T)$ is defined by $v(x_p(t),t)=\dot{x}_p(t).$

For a function $f(t) : \Gamma(t) \to \mbR$, $t \in [0,T)$, the material time derivative in a point $(x,t)$ with $x \in \Gamma(t)$ is defined by
\begin{equation}
 \mtd{v} f (x,t) := \frac{d}{dt} f(x_p(t),t) = \frac{\partial \tilde{f}}{\partial t}(x,t) + v(x,t) \cdot \nabla \tilde{f}(x,t)
\end{equation}
where $x = x_p(t)$ for a material particle $p$ located at $x_p(t)$ at time $t$. Note that although the expressions on 
the right hand side require 
a smooth extension $\tilde{f}$ of $f$ to a neighbourhood of the evolving surface in order to be well-defined, the resulting 
derivatives are independent of the extension. The tangential or surface gradient is defined as the projection of the standard 
derivative onto the tangent plane of the surface so that
\[
\nabla_{\Gamma(t)} f(x,t) := \nabla \tilde{f}(x,t) - (\nabla \tilde{f}(x,t) \cdot \nu(x,t)) \nu(x,t).
\]
Analogously for the surface divergence of a vector valued surface field. We can define the Laplace-Beltrami operator as 
$\Delta_{\Gamma(t)} f := \nabla_{\Gamma(t)} \cdot \nabla_{\Gamma(t)} f$. For calculus on surfaces, integration by parts formulas, 
and transport identities involving the material time derivative we refer to \cite{DziEll13a}. 

\medskip

As a \emph{model problem} 
we consider a conserved field on an evolving surface subject to Fickian diffusion. 
The advection diffusion equation is, see \cite{DziEll07}, %
\begin{subequations}
 \begin{align}
 \label{eq:2}
 \mtd{v} u + u \nabla_{\Gamma(t)} \cdot v - D \Delta_{\Gamma(t)} u & = 0 && \mbox{ on } \Gamma(t), \, t \in (0,T), \\
 u( \cdot, 0 ) & = u_0 && \mbox{ on } \Gamma(0).
 \end{align}
\end{subequations}
with initial value $u_0 : \Gamma(0) \to \mbR$ and diffusion parameter $D>0$.
%
%
We refer to \cite{AlpEllSti15} for a weak well-posedness analysis of such equations including suitable time-dependent Sobolev spaces.

\subsection{Numerical scheme}

Recall that the scanning process creates a representation of the cell boundary at each time in terms of a closed evolving triangulated hypersurface. We assume that by increasing the computational effort, i.e., refining the triangulation and the time step size, the resulting triangulations better approximate a closed evolving smooth hypersurface. For simplicity, we use linear surface finite elements and discretise in time by backward Euler. See \cite{Dem2009} for higher order approximations in space and \cite{LubManVen2013} for a discussion of time stepping schemes.

The triangulated surfaces given at times $t^{(m)}$, $m = 0, \dots, M$ with $t^{(M)} = T$, are denoted by $\Gamma_h^{(m)}$. 
Here, $h$ stands for the diameter of the largest triangle of the $\Gamma_h^{(m)}$. Let $S_h^{(m)}$ denote the space of piecewise linear finite element 
functions over $\Gamma_h^{(m)}$ with nodal basis $\{ \chi_j^{(m)} \}_{j=1}^J$. Time steps are denoted by $\tau^{(m)} = t^{(m)} - t^{(m-1)}$, $m = 1, \dots, M$. 
The vertex velocity of the triangulation is denoted by $w_h^{(m)} \in (S_h^{(m)})^3$ and it is defined by choosing the vertex values to be 
$w_h^{(m)}(x_j^{(m)}) = (x_j^{(m)} - x_j^{(m-1)}) / \tau^{(m)}$ where the $x_j^{(m)}$ are the vertex positions of $\Gamma_h^{(m)}$.

Given an approximation $u_h^{(0)}$ to the initial values $u_0$, the finite element discretisation of \eqref{eq:2} is to subsequently 
find $u_h^{(m)} \in S_h^{(m)}$, $m=1, \dots, M$, of the form $u_h^{(m)} = \sum_{j=1}^J u_j^{(m)} \chi_j^{(m)}$ as the solution to
\begin{multline} \label{eq:discrete}
\frac{1}{\tau^{(m)}} \left( \int_{\Gamma_h^{(m)}} u_h^{(m)} \chi_j^{(m)} \, \mathrm{d} \sigma_h - \int_{\Gamma_h^{(m-1)}} u_h^{(m-1)} \chi_j^{(m-1)} \, \mathrm{d} \sigma_h \right) \\
+ D \int_{\Gamma_h^{(m)}} \nabla_{\Gamma_h^{(m)}} u_h^{(m)} \cdot \nabla_{\Gamma_h^{(m)}} \chi_j^m \, \mathrm{d} \sigma_h + b_{adv} (u_h^{(m)}, \chi_j^{(m)}) + b_{sld} (u_h^{(m)}, \chi_j^{(m)}) = 0
\end{multline}
for $j = 1,2, \ldots, J$. Here, $b_{adv}$ and $b_{sld}$ are two terms which account for the material velocity $v$ in \eqref{eq:2} and for streamline diffusion, respectively. 

The advection term is
\[
 b_{adv} (u_h^{(m)}, \chi_j^{(m)}) = \int_{\Gamma_h^{(m)}} u_h^{(m)} \, \big{(} w_h^{(m)} -v_h^{(m)} \big{)} \cdot \nabla_{\Gamma_h^{(m)}} \chi_j \mathrm{d} \sigma_h.
\]
Here, $v_h^m$ denotes the discrete material velocity. 

In ALE methods a non-physical arbitrary velocity for the vertices may be used to ensure a better quality of the evolving mesh \cite{EllSty12,EllVen15}. 
In our setting the vertex velocity is dependent on the computational tools used to obtain the parametrisation and hence may not be physical. In our case where the material velocity is 
purely in the normal direction we set $v_h^{(m)} = w_h^{(m)} \cdot \nu_h^{(m)} \nu_h^{(m)}$ on each triangle where $\nu_h^{(m)}$ is the (piecewise constant) unit normal of the triangle. 

In some applications the tangential vertex velocity is quite strong in comparison with transport by diffusion, i.e., the problem is advection dominated which can 
lead to poor approximations on coarse grids. A standard way to deal with this issue is to add a streamline diffusion term, \cite{HugBro79}, which is of the form 
\begin{equation*}
b_{sld} (u_h^{(m)}, \chi_j^{(m)}) = D g(h) \int_{\Gamma_h^{(m)}} 
( w_h^{(m)} \cdot \nabla_{\Gamma_h^{(m)}} u_h^{(m)} ) ( w_h^{(m)} \cdot \nabla_{\Gamma_h^{(m)}} \chi_j^{(m)} ) \, \mathrm{d} \sigma_h
\end{equation*}
where we took $g(h) = h^2$ in the computations below. In each time step $m$, \eqref{eq:discrete} yields a linear system of 
equations for coefficients $u_j^{(m)}$ which can be solved using a biconjugate gradient stabilised method.


%% file: source/application_submitted.tex
\section{Applications and simulations}

In this section we illustrate the use of evolving surface finite element methodology on triangular meshes obtained from image data 
of neutrophil cells. These  were labelled with cell mask orange dye to stain the plasma membrane. 
The neutrophils are quiescent initially and then stimulated by application of fMLP which prompts cell spreading and movement \cite{HawSteSuiWil11}. 
Cell movements are imaged by spinning disk confocal microscopy \cite{GraRietZim05}. The acquisition speed is about 80 ms per slice (4 secs/stack). 

The data comprises of 166 surface triangulations associated with a subset of 214 successive times where each triangular mesh has $20 480$ triangles with $10 242$ vertices. In the computations triangulations needed for intermediate times were obtained by interpolation of the vertex positions. 


\subsection{Simulation of FRAP experiments}
\label{sec:frapsim}

Fluorescence recovery after photobleaching (FRAP) is a standard technique to assess the mobility of fluorescently labelled molecules in 
cells; for example, cell surface receptors which bind extracellular chemoattractants such as fMLP. Within a small region these molecules 
are irreversibly photo-bleached by a high intensity laser impulse. This is then followed by a fluorescence recovery phase consisting 
of the motion of non-affected molecules into the bleached region. Monitoring this gives insight into the type of motion and eventually 
enables the quantification of kinetic parameters. 
We refer to \cite{ElzJan2009} for an overview of what type of information can be obtained from FRAP. 

Previous studies have considered models of this recovery process on a flat stationary surface. For example, the mathematical 
modelling in \cite{HagRoeCruBar2005} uses semi-analytical and computational finite difference methods which do not account for the complex 
shape of the membrane and its evolution. Here, we show that using an accurate representation of the geometry affects the recovery of the 
concentration, denoted by $c$, in the bleached region. 
Our model for FRAP consists of the diffusion equation of the form \eqref{eq:2}. In particular, we neglect effects such as non-instantaneous 
bleaching and cytoplasmic attachment/detachment kinetics which sometimes are of importance \cite{HagRoeCruBar2005}.
We assume a model in which there is no tangential, advective material transport so that the material velocity is taken in the normal direction.
Hence we have $$v_h^{(m)} = w_h^{(m)} \cdot \nu_h^{(m)} \nu_h^{(m)}.$$

For the diffusion coefficient we take the values 
$D = 5 \cdot 10^{-10} \, \mathrm{cm}^2 / \mathrm{sec}$ which is a typical value for G-protein coupled receptors such as fMLP. 
The initial photobleaching is modelled by a homogeneous constant value $C_0$ except in a  \emph{region of interest} (ROI) $\Gamma_b(0)$,
defined by the  intersection of $\Gamma(0)$ with an idealised
ball of radius $r$ centred in $m$. 
Here $r$ is small with respect to the diameter of the cell enclosed by $\Gamma(0)$. 
We remark that realistic bleach profiles and side effects, such as photobleaching by light scattering, due to choosing this (ROI) are neglected 
and refer to \cite{MazCelVicKroDia07} for a discussion of these issues.
In our setting the  observation time is short and  we do not  advect the the bleached material region with the surface velocity.
Our region of observation of the concentration is thus  $\Gamma_b(t)$ defined by 
the intersection of $\Gamma(t)$ with the (ROI), the fixed   ball of radius $r$ centred at $m$. 
Our approach mimics that of typical experimental set ups where fluorescence recovery is sampled within a fixed (ROI).

We rescale the concentration so that $C_0 = 1$ and 
show concentration relative to $C_0$. Precisely, the initial condition is given by
\begin{equation*}
 c(0,x) = c_0(x) =
 \begin{cases}
 0 & \mbox{if } x \in \Gamma_b = \{ x \in \Gamma : | x - m | < r \} \\
 1 & \mbox{otherwise}.
 \end{cases}
\end{equation*}
In the computations, we take a piecewise linear interpolation of $c_0$ for initial condition. 
The time step is chosen to be $0.04 \, \mathrm{s}$ - i.e. 100th of the difference between successive  images.

We present three simulations on different surface evolutions using our single data set starting from frames 0, 25 and 77.
In each case we take the (ROI)  to be a sphere of radius a quarter the radius of the cell centred at $(0.25, 0.25, 0.25)$. The  areas of bleached regions are, respectively, $8\%$, $5\%$ and $3\%$ of the total cell areas.
In the calculations  $\Gamma_b(t)$ is approximated by 
the evolution of the  union of elements comprising the initial discrete (ROI).

\begin{figure}
 \centering
 \begin{subfigure}[b]{\textwidth}
 \includegraphics[width=0.24\textwidth]{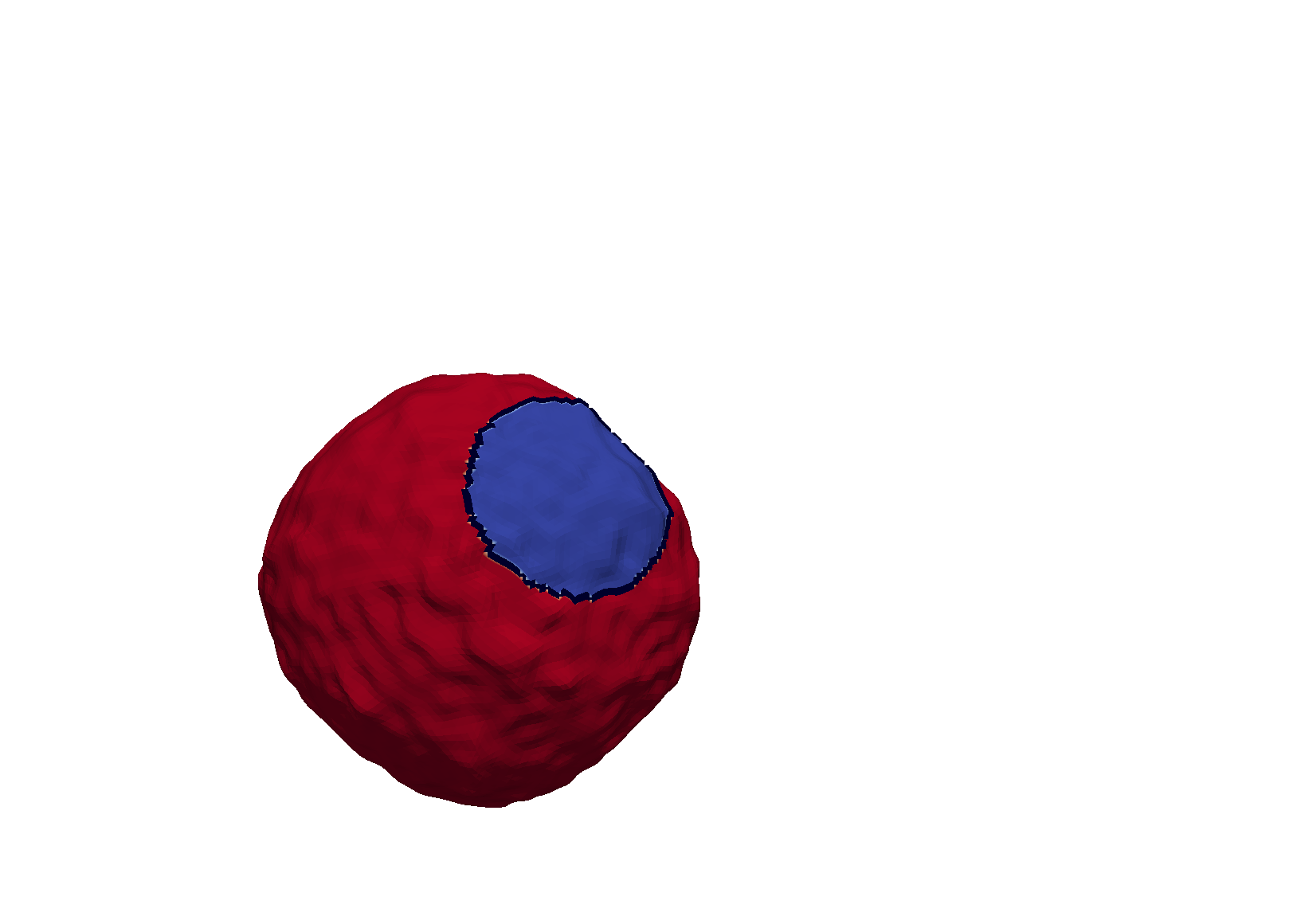}
 \includegraphics[width=0.24\textwidth]{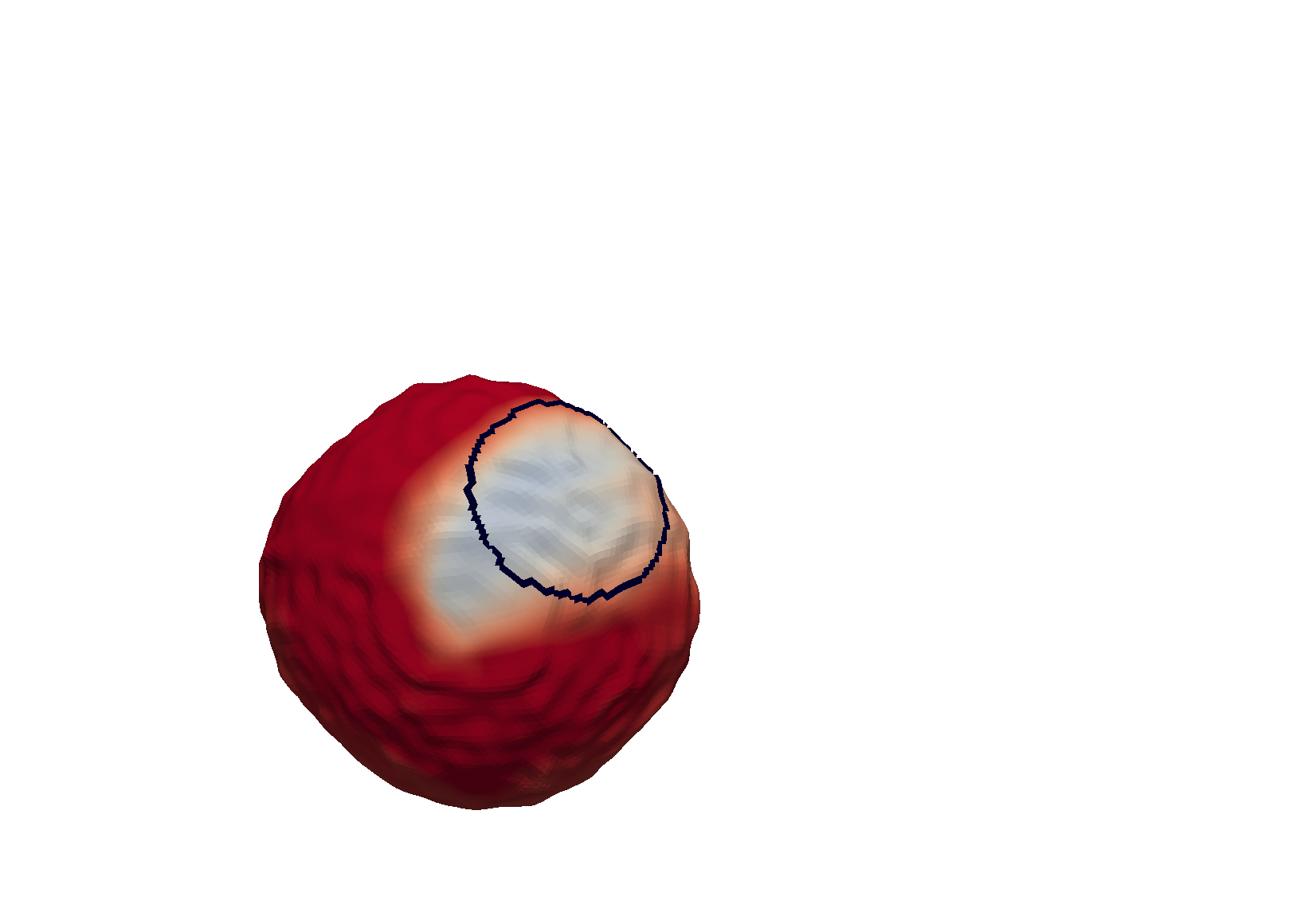}
 \includegraphics[width=0.24\textwidth]{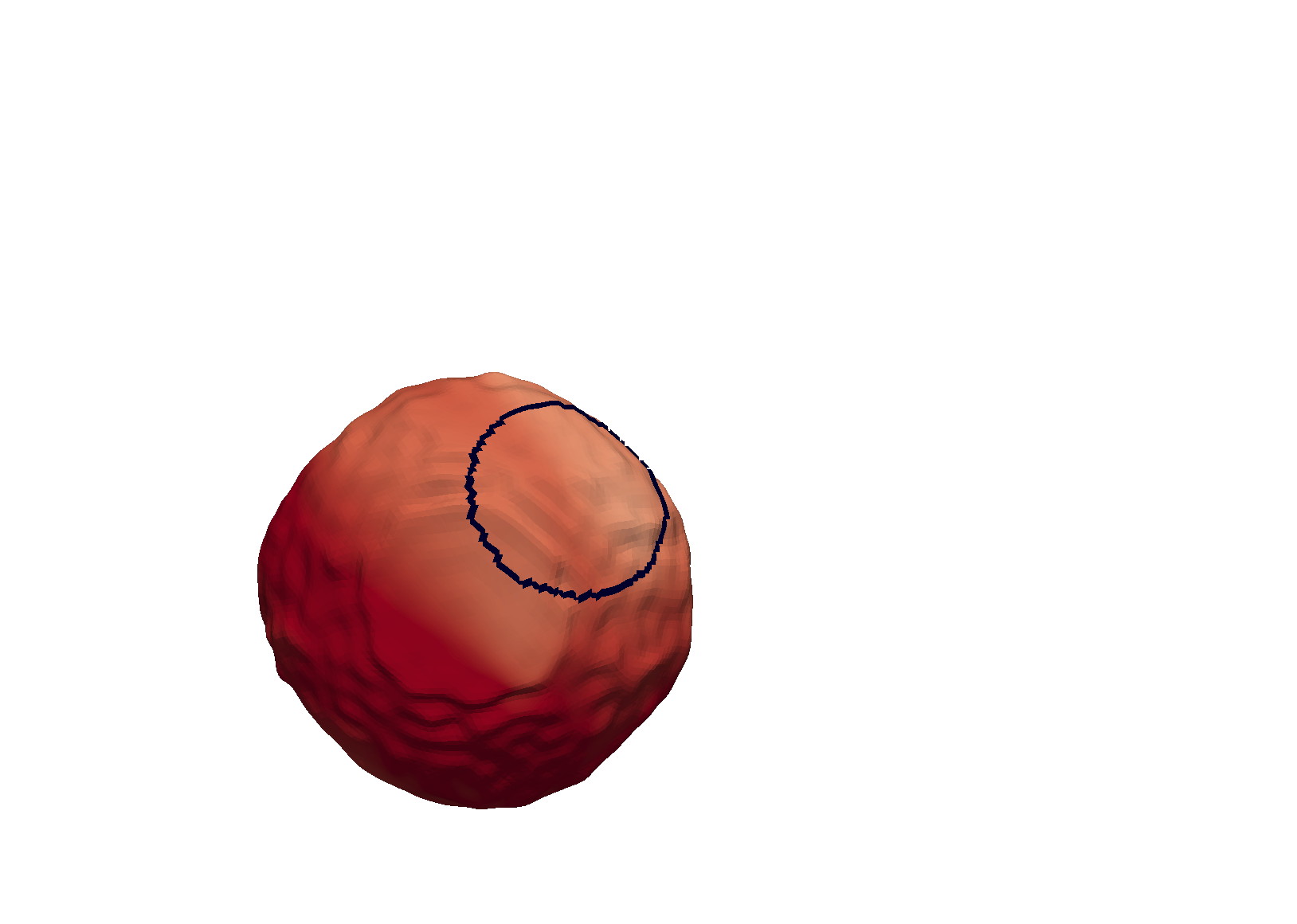}
 \includegraphics[width=0.24\textwidth]{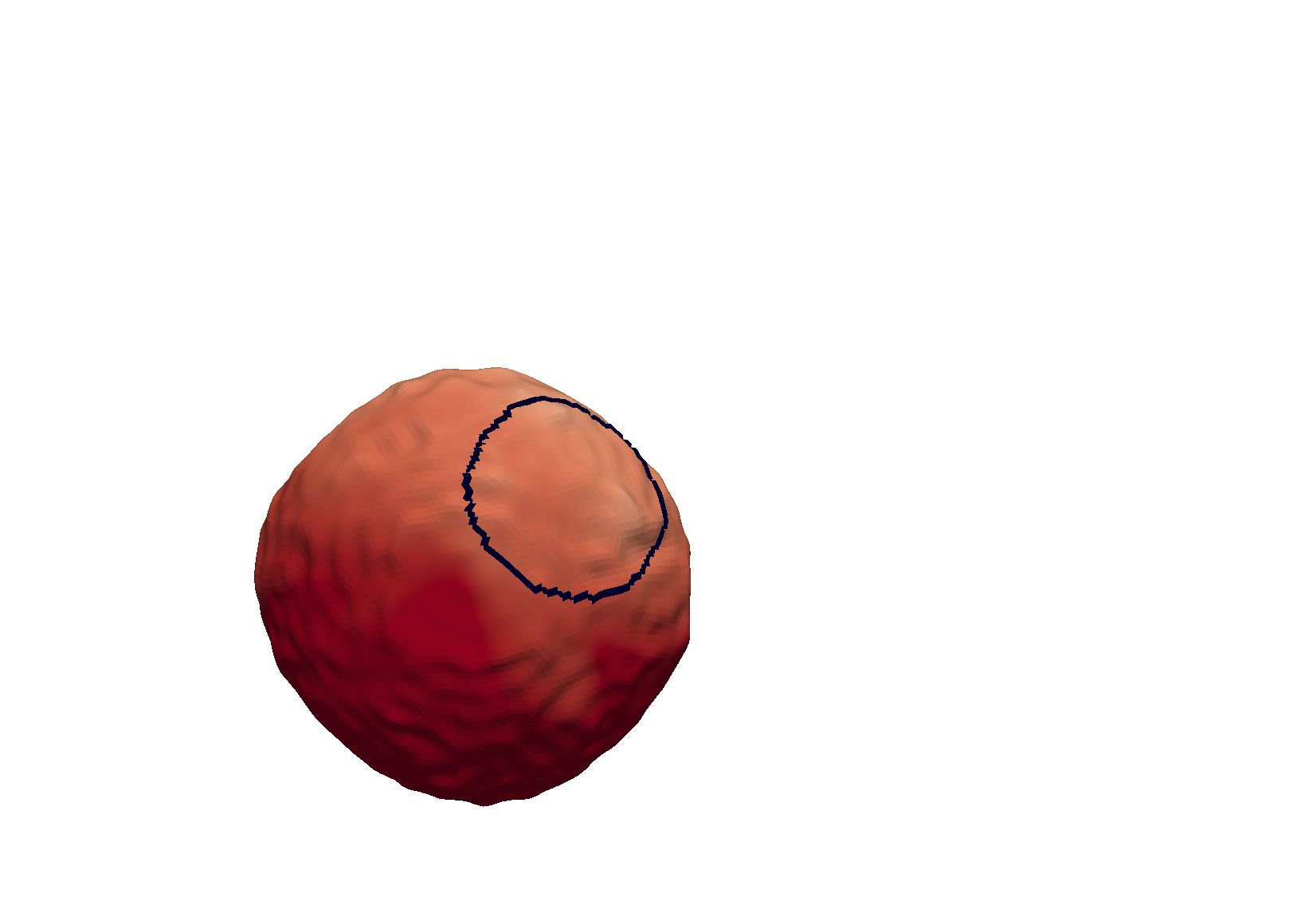}
 \caption{Starting at frame 0. Plots at times $t = 0, 4, 8, 12 \mathrm{s}$ from bleaching.}
 \label{fig:frapA-normal}
 \end{subfigure}

 \begin{subfigure}[b]{\textwidth}
 \includegraphics[width=0.24\textwidth]{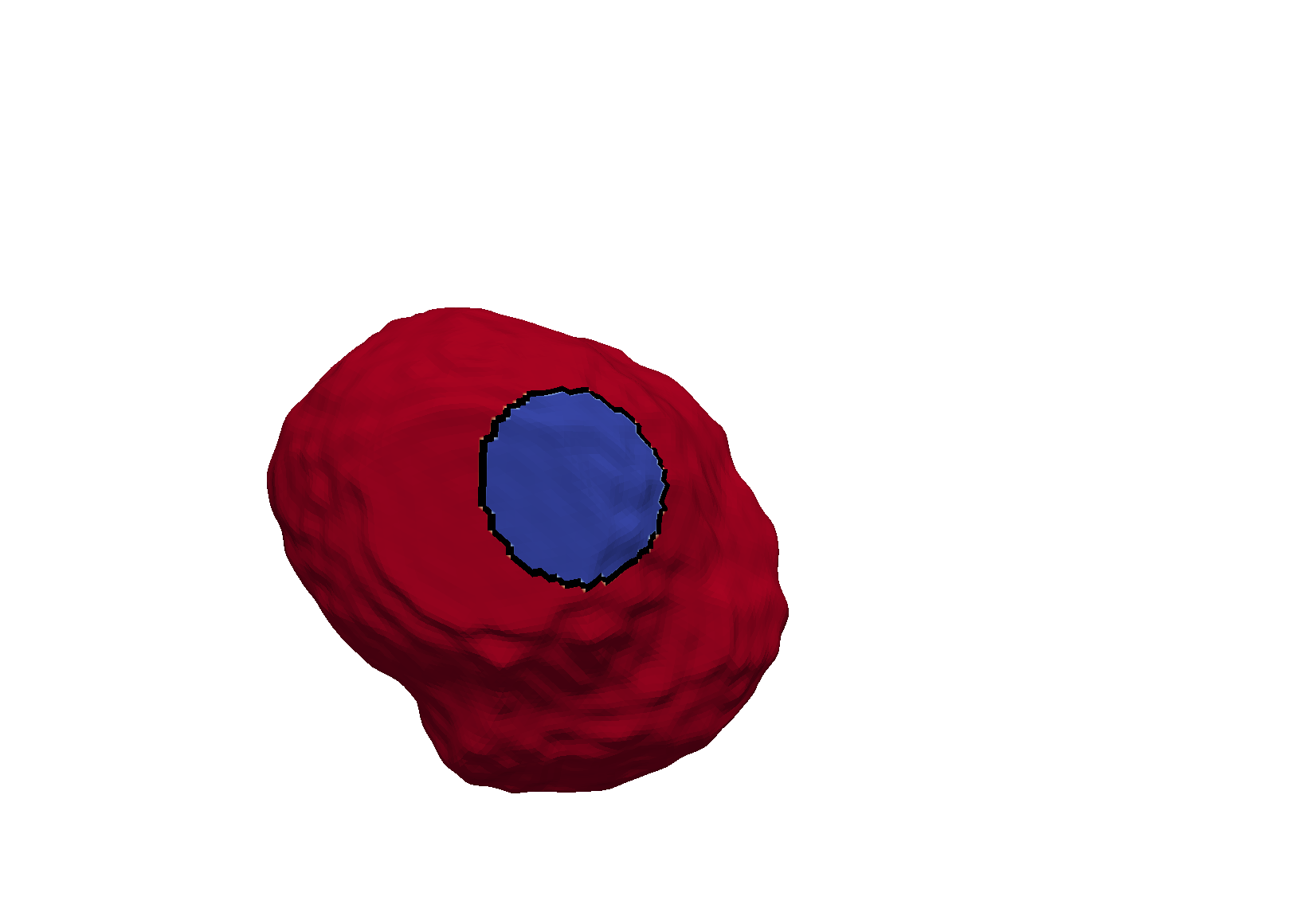}
 \includegraphics[width=0.24\textwidth]{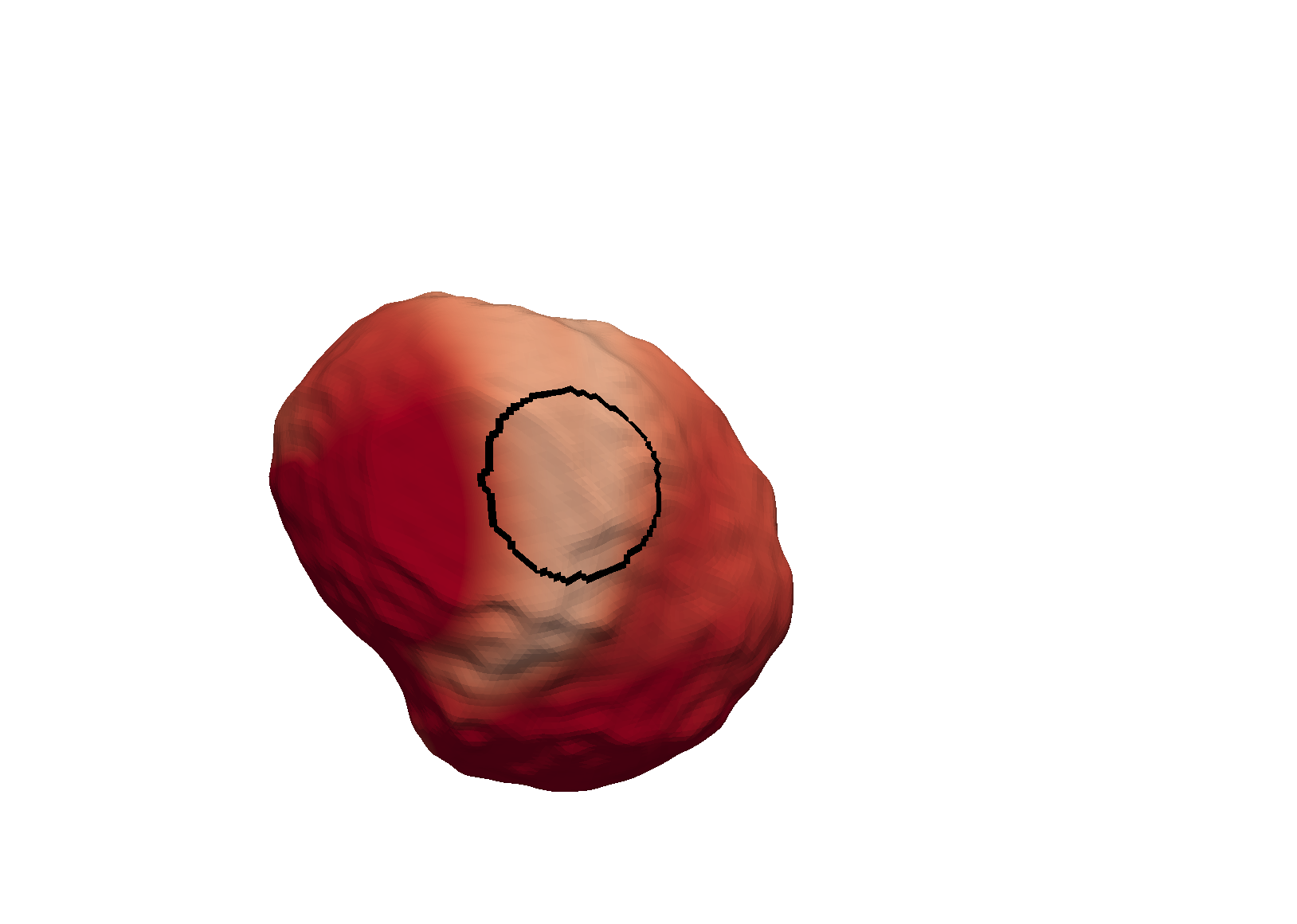}
 \includegraphics[width=0.24\textwidth]{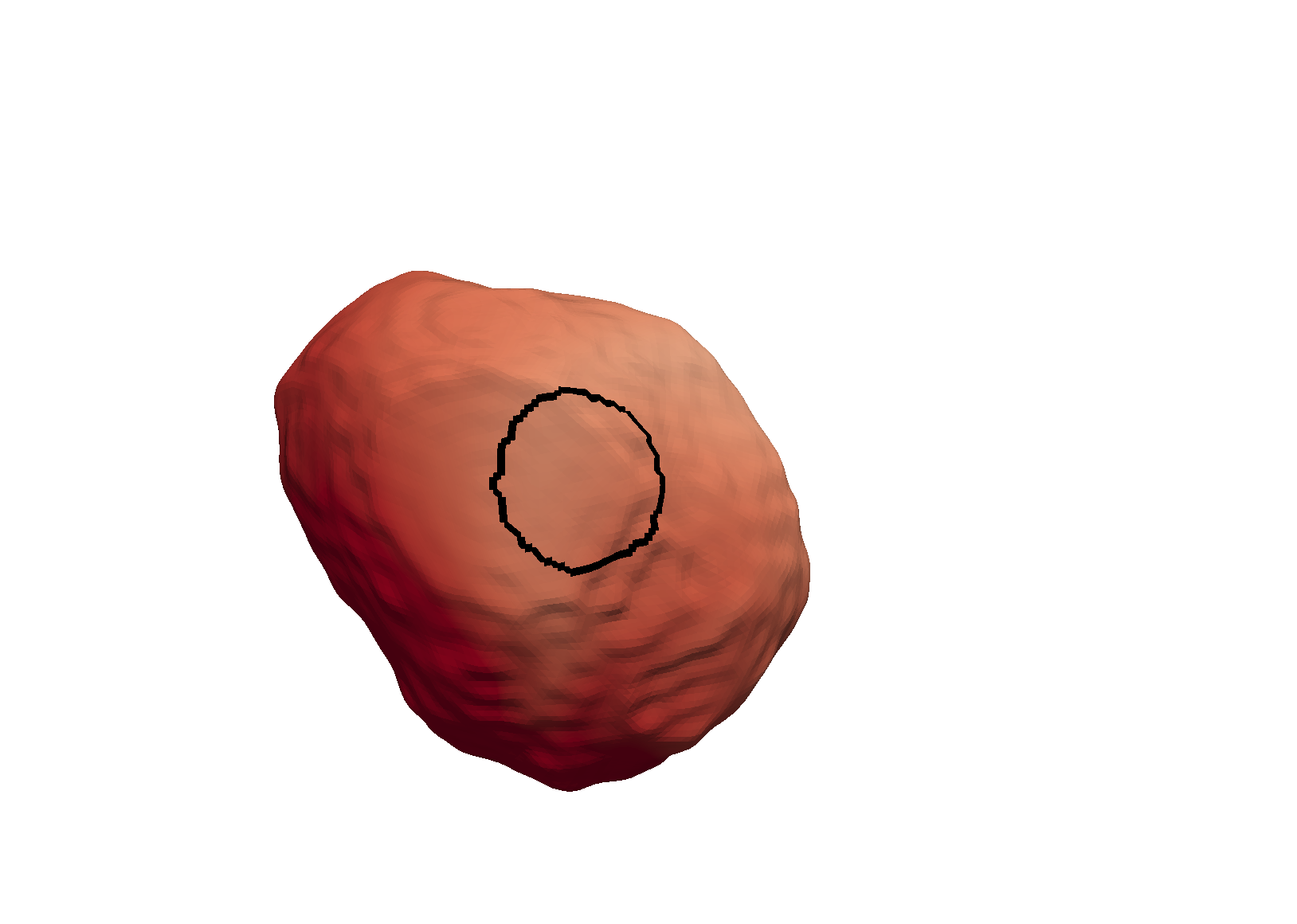}
 \includegraphics[width=0.24\textwidth]{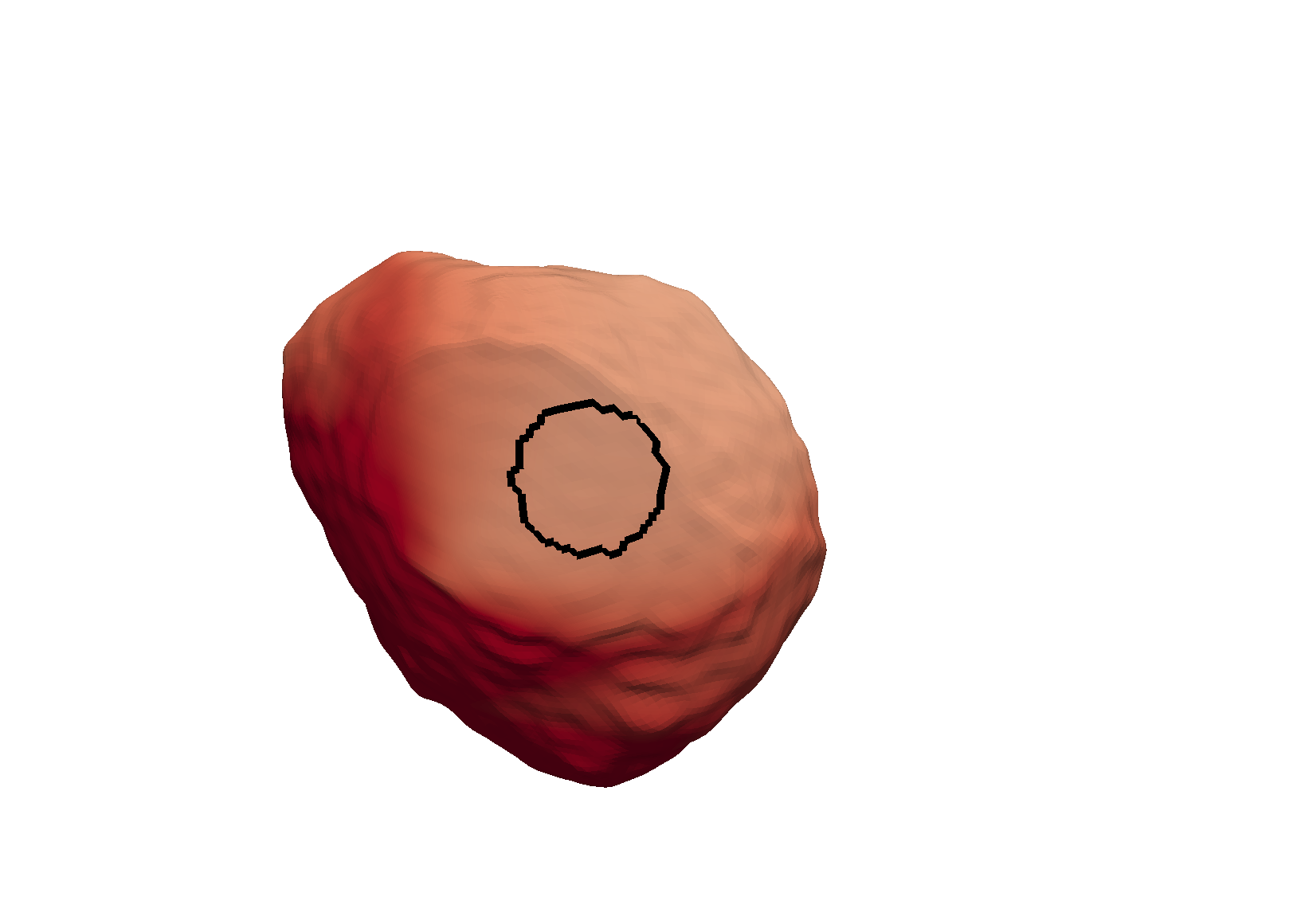}
 \caption{Starting at frame 25. Plots at times $t = 0, 4, 8, 12 \mathrm{s}$ from bleaching.}
 \label{fig:frapB-normal}
 \end{subfigure}

 \begin{subfigure}[b]{\textwidth}
 \includegraphics[width=0.24\textwidth]{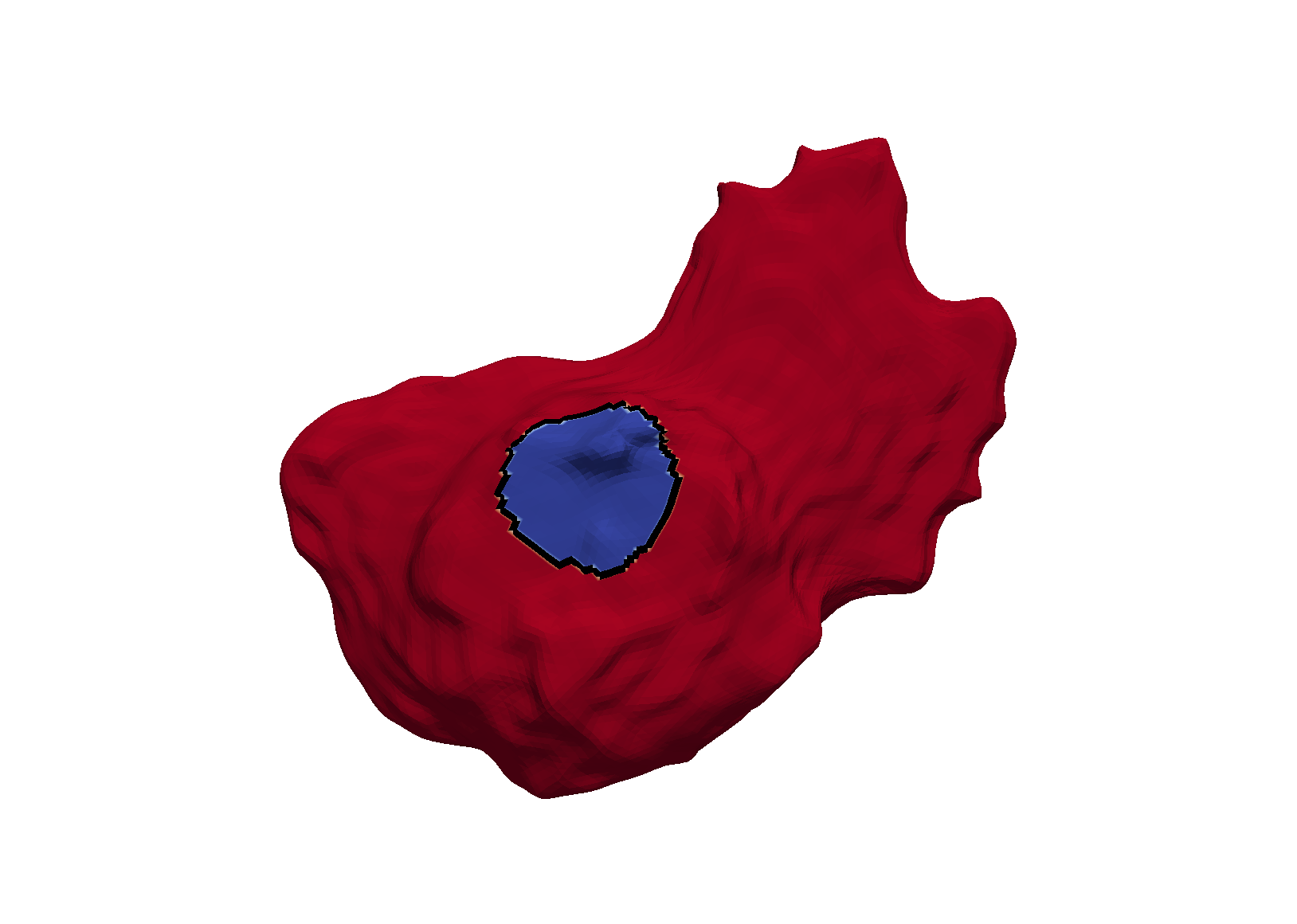}
 \includegraphics[width=0.24\textwidth]{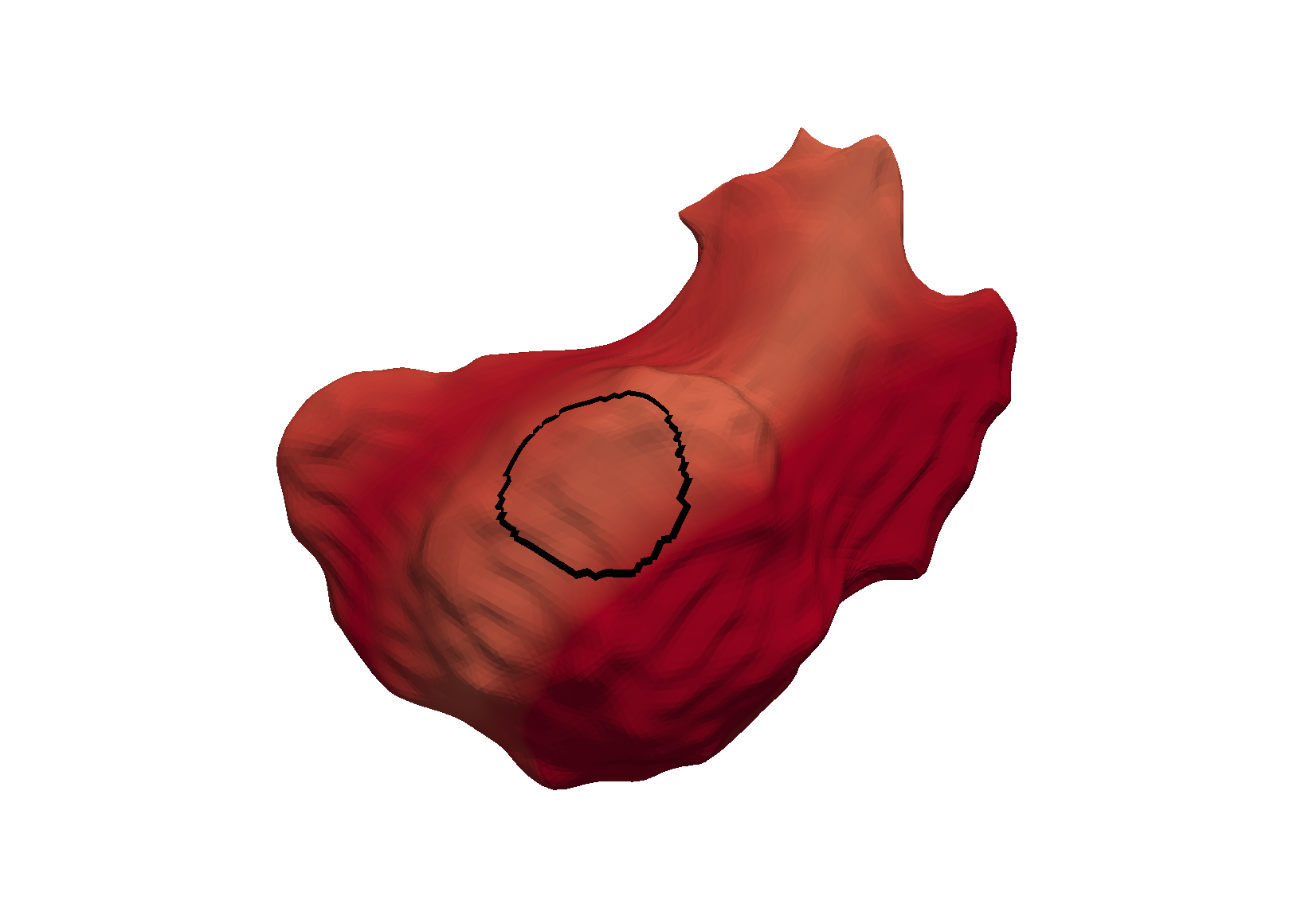}
 \includegraphics[width=0.24\textwidth]{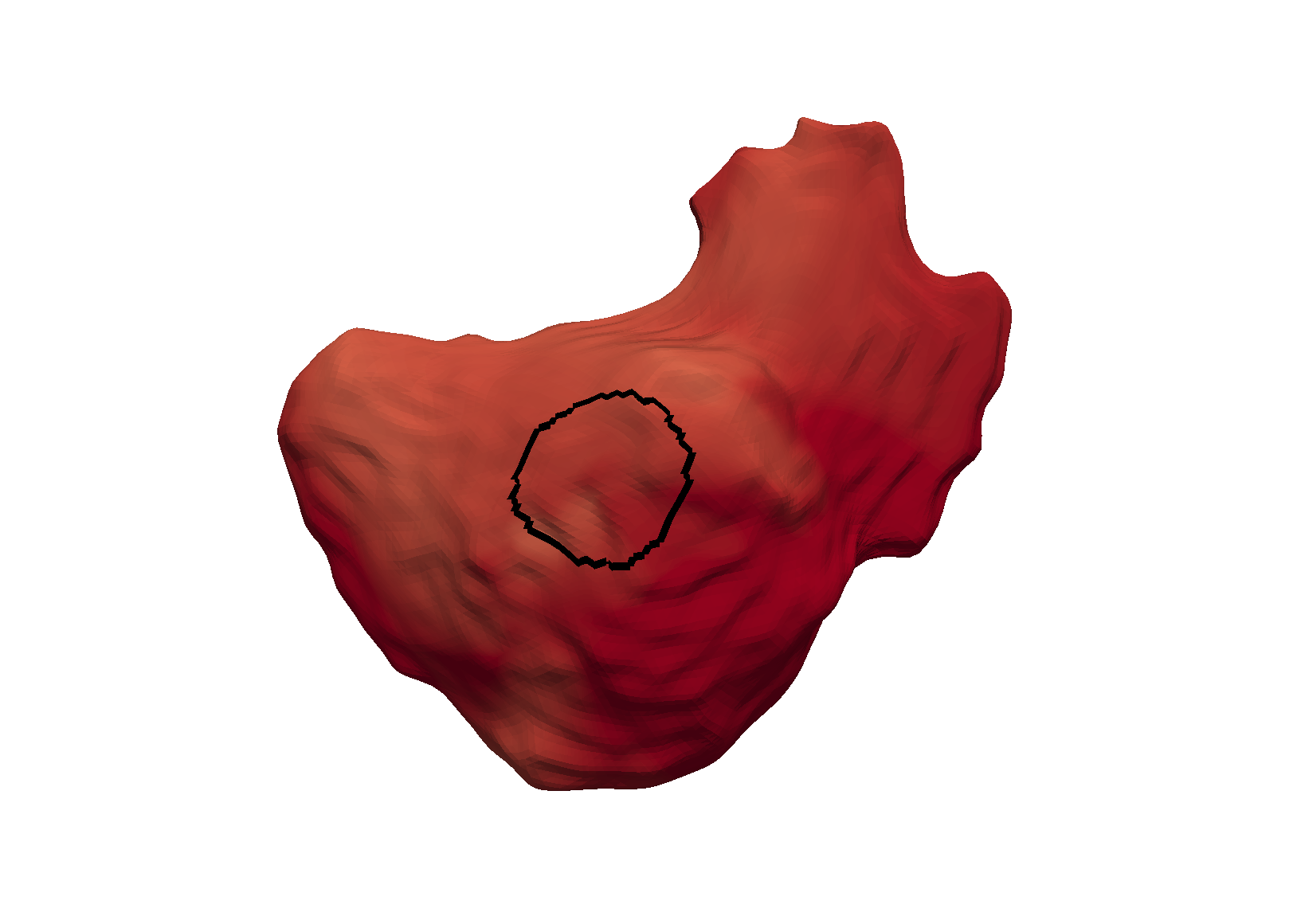}
 \includegraphics[width=0.24\textwidth]{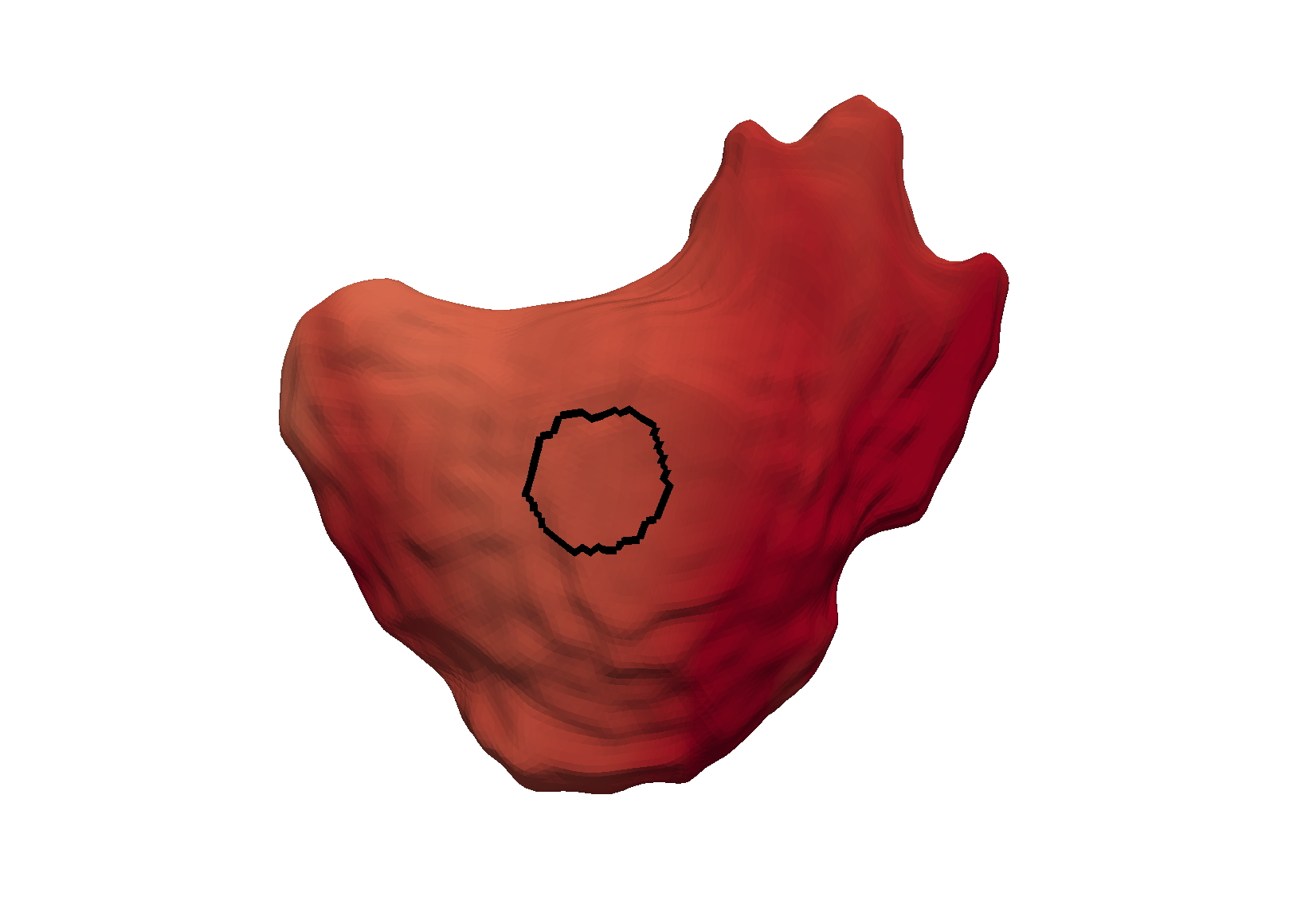}
 \caption{Starting at frame 77. Plots at times $t = 0, 4, 8, 12 \mathrm{s}$ from bleaching.}
 \label{fig:frapC-normal}
 \end{subfigure}

 \caption{Results of FRAP bleaching simulations. The colour scheme is from red (1 arb.) to blue (0 arb.). 
 A black contour shows the boundary of the (ROI) $\Gamma_b(t)$ in which we track the concentration over time.}
 \label{fig:frap-results-normal}
\end{figure}

Snapshots from each simulation are shown in Figure~\ref{fig:frap-results-normal}. In the first simulation (Figure~\ref{fig:frapA-normal}), the cell  surface 
undergoes very small changes and has no large changes in curvature. In the other two simulations (Figures~\ref{fig:frapB-normal} and \ref{fig:frapC-normal}),
the overall cell surfaces display large deviations in curvature. These changes clearly change the rate at which the concentration recovers in $\Gamma_b$.

\begin{figure}

 \begin{subfigure}[c]{0.45\textwidth}
 \includegraphics{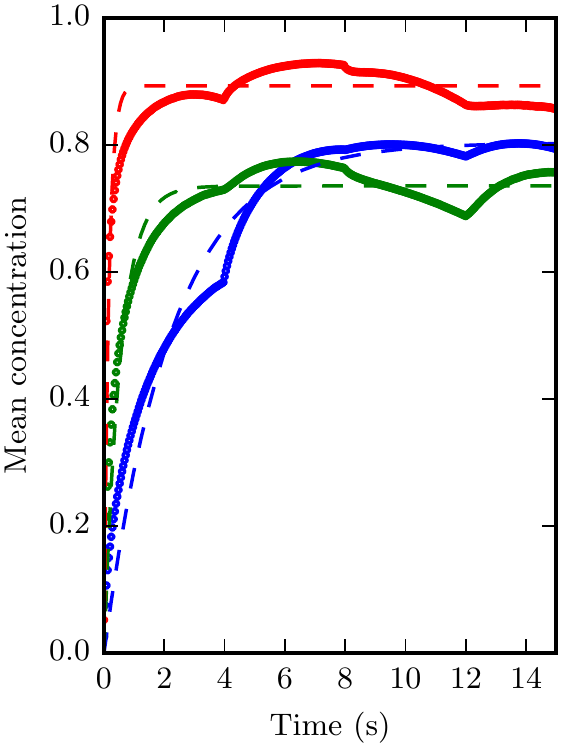}
 \end{subfigure}
 \begin{subfigure}[c]{0.54\textwidth}
 \begin{tabular}{c|c|c|c}
   Run & $A$ (arb.) & $B$ ($\mathrm{s}$) & $T_{1/2}$ ($\mathrm{s}$) \\
   \hline
   (A) & $0.80$ ($3.9 \cdot 10^{-3}$) & $2.3$ ($4.2 \cdot 10^{-2}$) & $1.6$ \\
   (B) & $0.73$ ($2.1 \cdot 10^{-3}$) & $0.54$ ($1.4 \cdot 10^{-2}$) & $0.38$\\
   (C) & $0.89$ ($2.1 \cdot 10^{-3}$) & $0.16$ ($6.4 \cdot 10^{-3}$) & $0.11$ \\
 \end{tabular}
 \end{subfigure}

 \caption{Curve fitting the recovery times for FRAP simulations.
   (Left) The recovery curves with: circles for observations;
   solid line for fitted curves;
   blue, green and red colour for geometries (A), (B) and (C), respectively.
   (Right) Fitted values with one standard deviation of the mean shown in brackets.
   We also report $T_{1/2}$ the time from the bleach to the timepoint where
   the fluorescence intensity reaches the half of the final recovered
   intensity. See Section \ref{sec:frapsim} for details on the simulation parameters.
   }
 \label{fig:recovery}
\end{figure}

In order to quantify the effects of the change in geometry we compute 
$\frac{1}{ | \Gamma_b (t) | }\int_{\Gamma_b (t)} c(t) \, \mathrm{d} \sigma$ at 
each time point. This quantity is then fitted to  a recovery curve $A ( 1 - \exp( -t / B ) )$ to the first $12 \, \mathrm{s}$ of simulation time (i.e. using four images).
We use this  short time because of the relatively small evolution of the overall cell boundary. The  underlying 
cell geometry is rather different in each case. It follows that  $A$ is the asymptotic mean value of the concentration in the evolved bleached region 
and  $B$ is then the recovery timescale. 
This is performed using the \verb!scipy.optimize! function \verb!curve_fit!. 
The results are shown in Figure~\ref{fig:recovery}. The range of resulting recovery rates show a large variation with differences of an order of magnitude.

\subsection{Pattern formation on evolving surfaces}
\label{sec:pattern}

Systems of reaction-diffusion equations such as Turing systems which involve local autocatalysis and inhibition on 
a longer range are well established tools to explain biological pattern formation. Recently, the behaviour of such systems on 
evolving domains has been increasingly attracting interest \cite{BerSbaKou10,BarEllMad11,VenLakMad2012}. 
Without any specific application in mind but to display the capability of our computational methodology we simulate an 
activator -- depleted substrate system taken from \cite{GieMei1972,BarEllMad11} for stripes or spots (in dependence of the parameters). In such systems, the activator, denoted by $u$, grows by consuming a substrate, denoted by $w$, which thanks to higher diffusivity acts as 
the long-range inhibitor. The partial differential equations on a moving surface read
\begin{equation} \label{eq:patternRDS}
 \begin{split}
 \mtd{v} u + u \nabla_{\Gamma(t)} \cdot v
 & = D_u \Delta_{\Gamma(t)} u + \gamma( a - u + u^2 w ) \\
 \mtd{v} w + w \nabla_{\Gamma(t)} \cdot v
 & = D_w \Delta_{\Gamma(t)} w + \gamma( b - u^2 w ).
 \end{split}
\end{equation}
We use the parameters
\begin{align*}
 D_u = 1, \quad
 D_w = 10, \quad
 \gamma = 200, \quad
 a = 0.1, \quad
 b = 0.9.
\end{align*}
For the initial values we consider $10\%$ random perturbations of the steady-state solution $(1,0.9)$.

\begin{figure}[p]
 \centering


\begin{subfigure}[b]{\textwidth}
\includegraphics[width=0.22\textwidth,clip=true,trim=15cm 10cm 15cm 10cm]{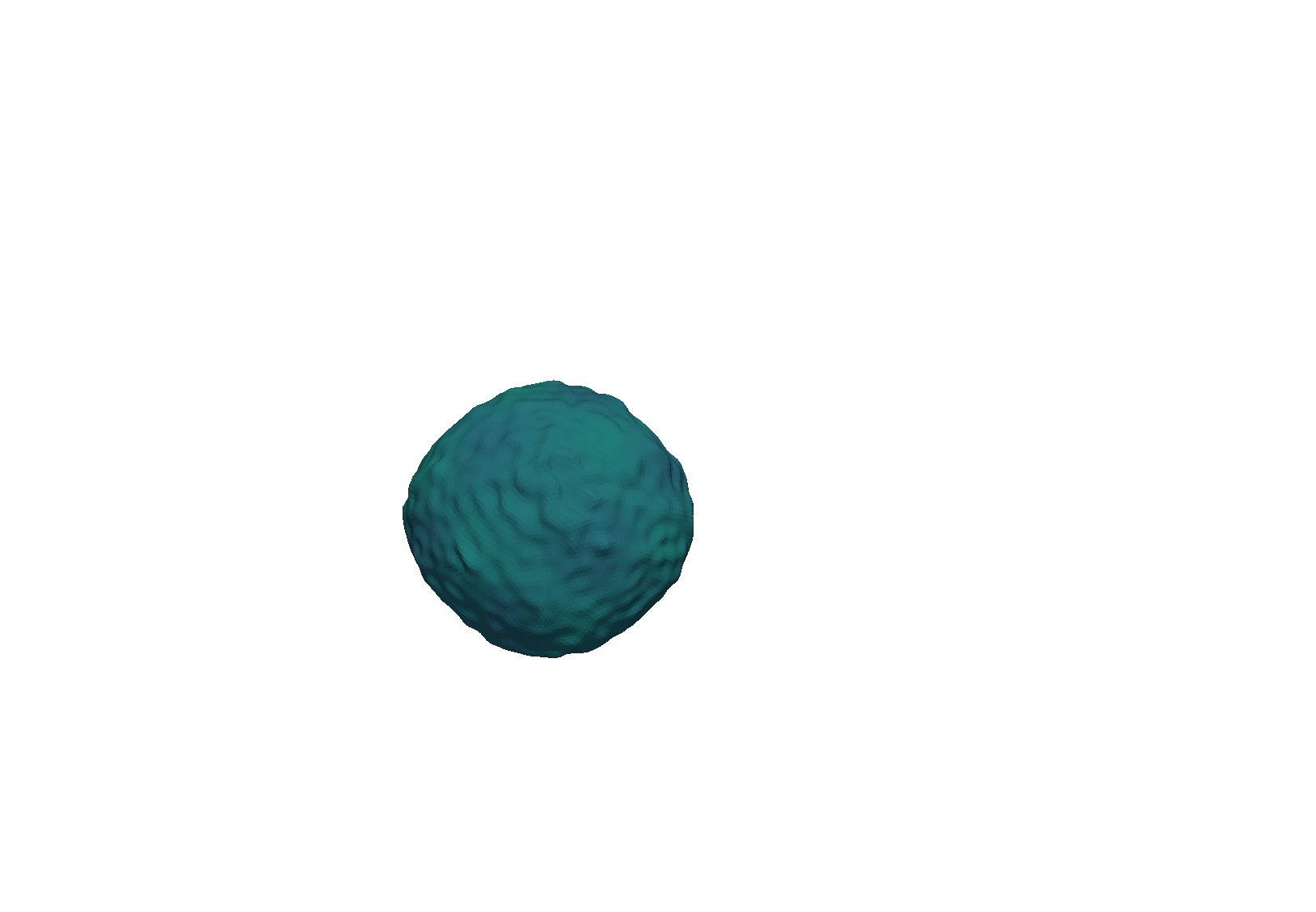}
\includegraphics[width=0.22\textwidth,clip=true,trim=15cm 10cm 15cm 10cm]{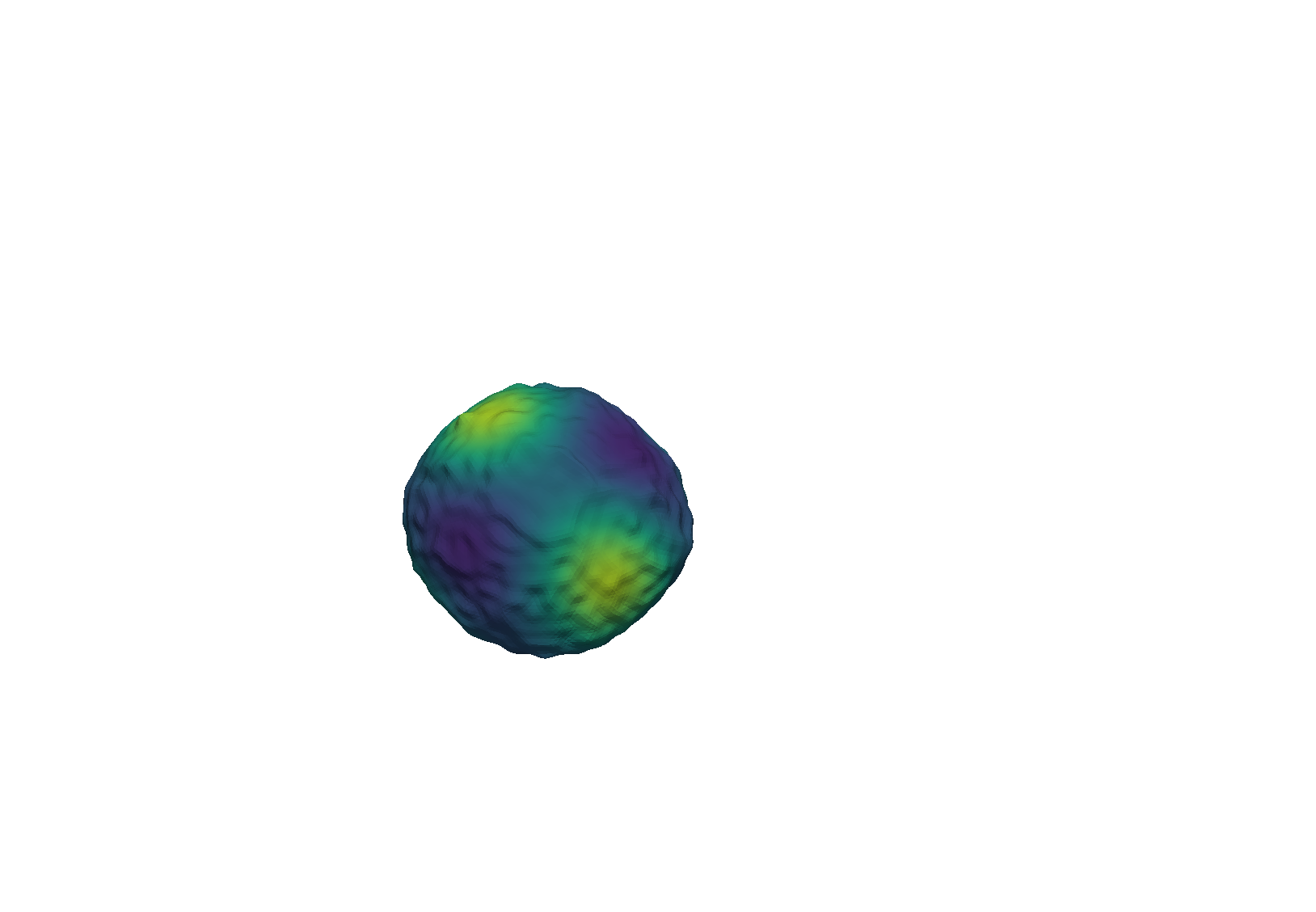}
\includegraphics[width=0.22\textwidth,clip=true,trim=15cm 10cm 15cm 10cm]{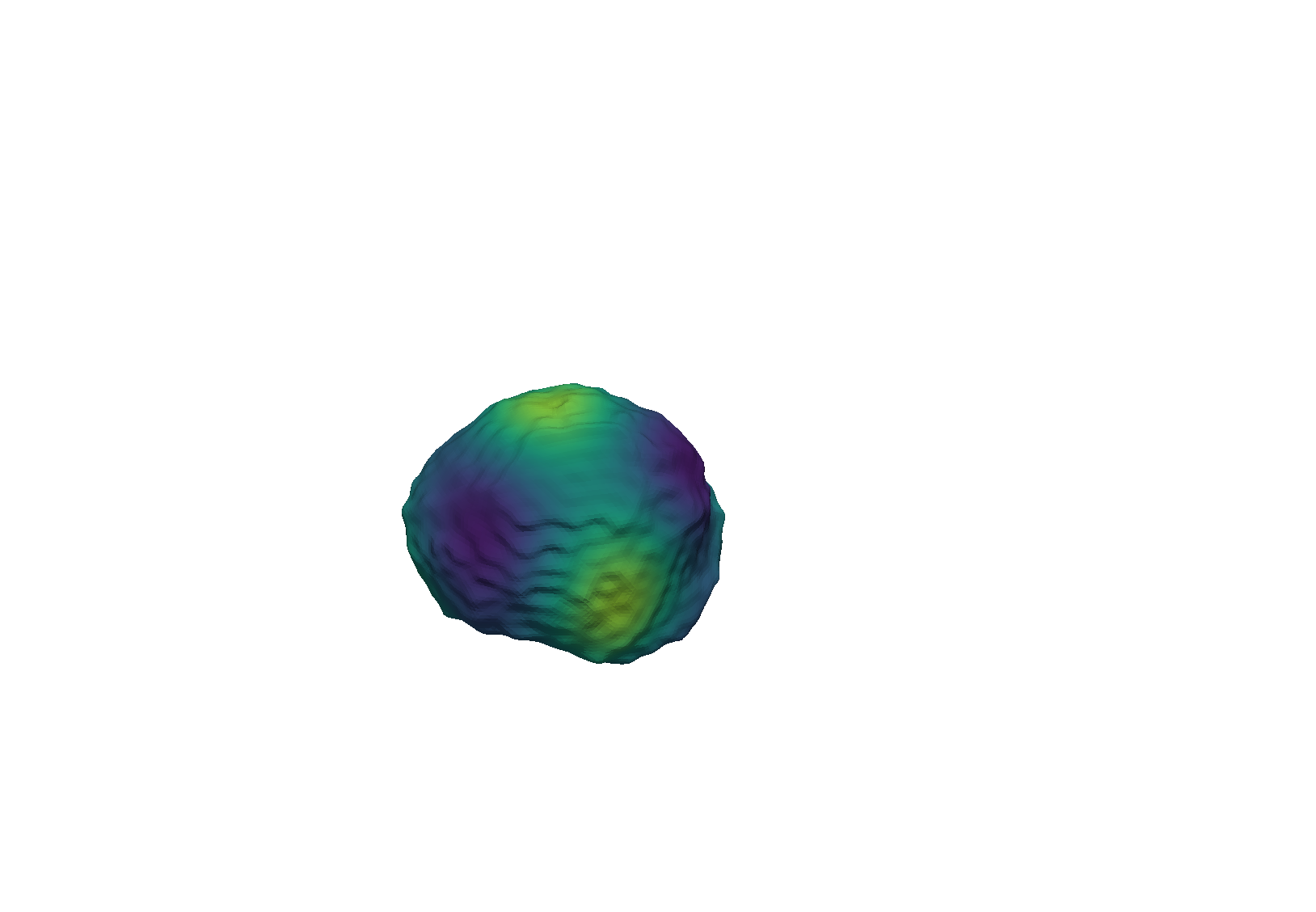}
\includegraphics[width=0.22\textwidth,clip=true,trim=15cm 10cm 15cm 10cm]{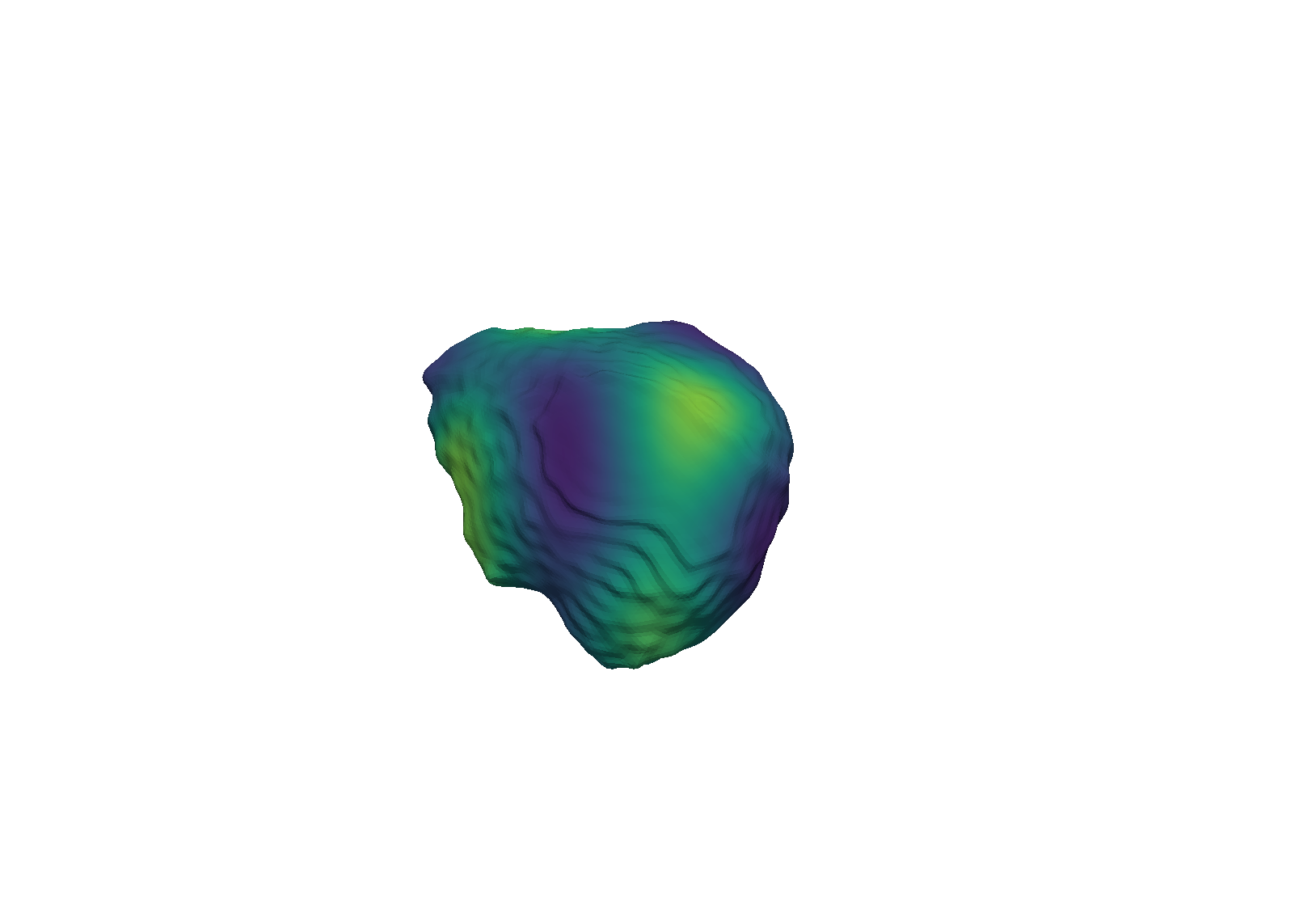}
\includegraphics[width=0.22\textwidth,clip=true,trim=15cm 10cm 15cm 10cm]{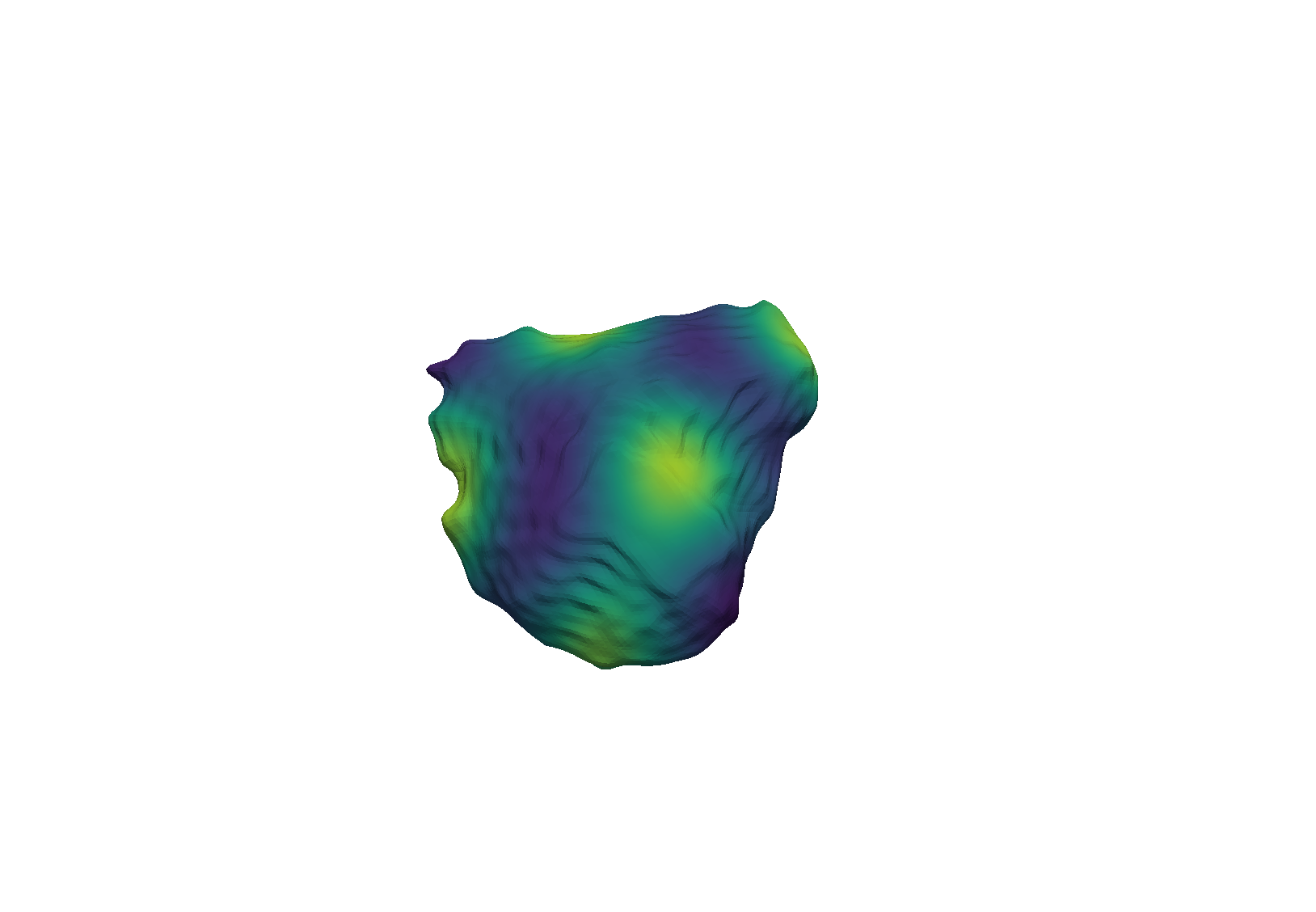}
\includegraphics[width=0.22\textwidth,clip=true,trim=15cm 10cm 15cm 10cm]{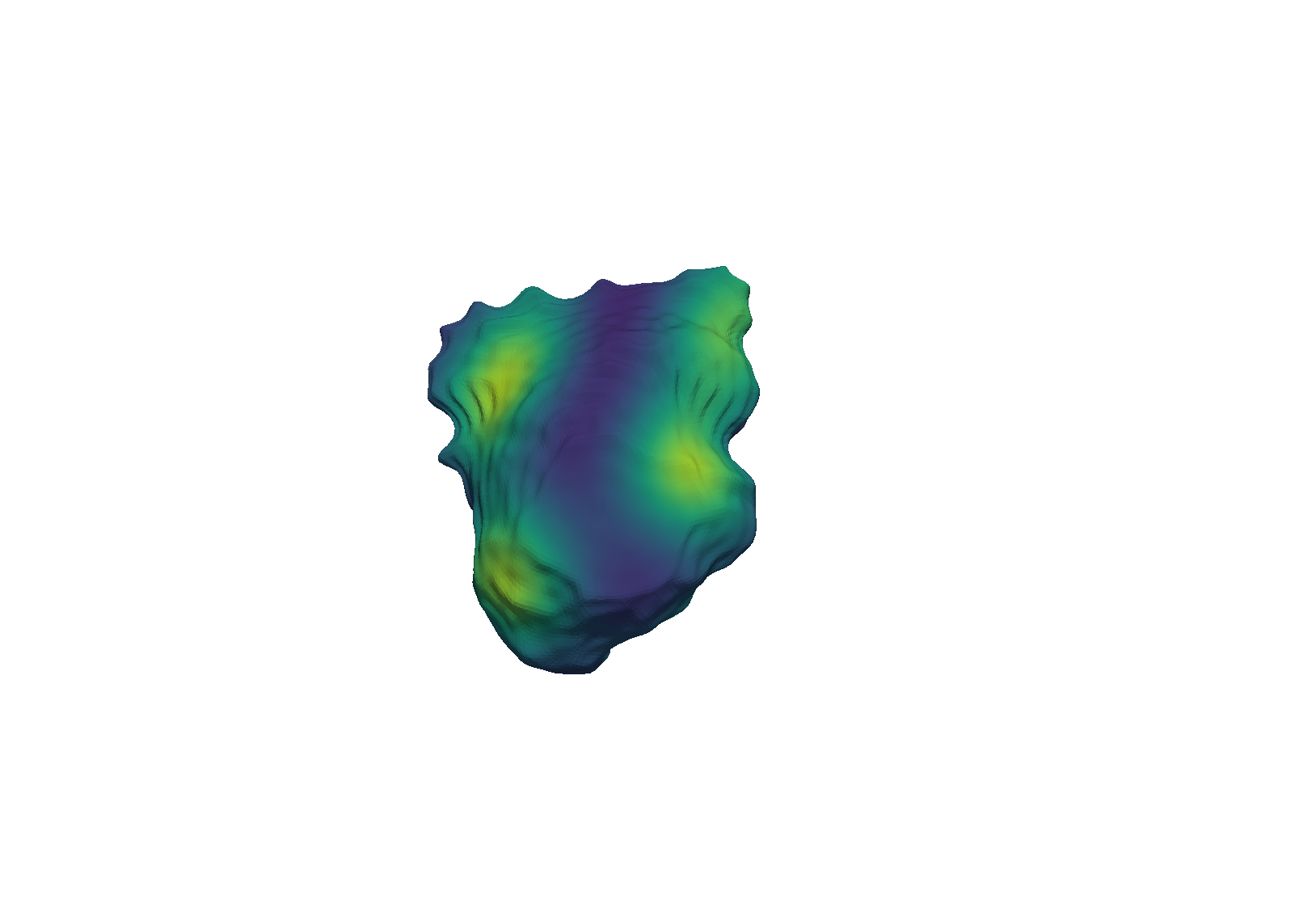}
\includegraphics[width=0.22\textwidth,clip=true,trim=15cm 10cm 15cm 10cm]{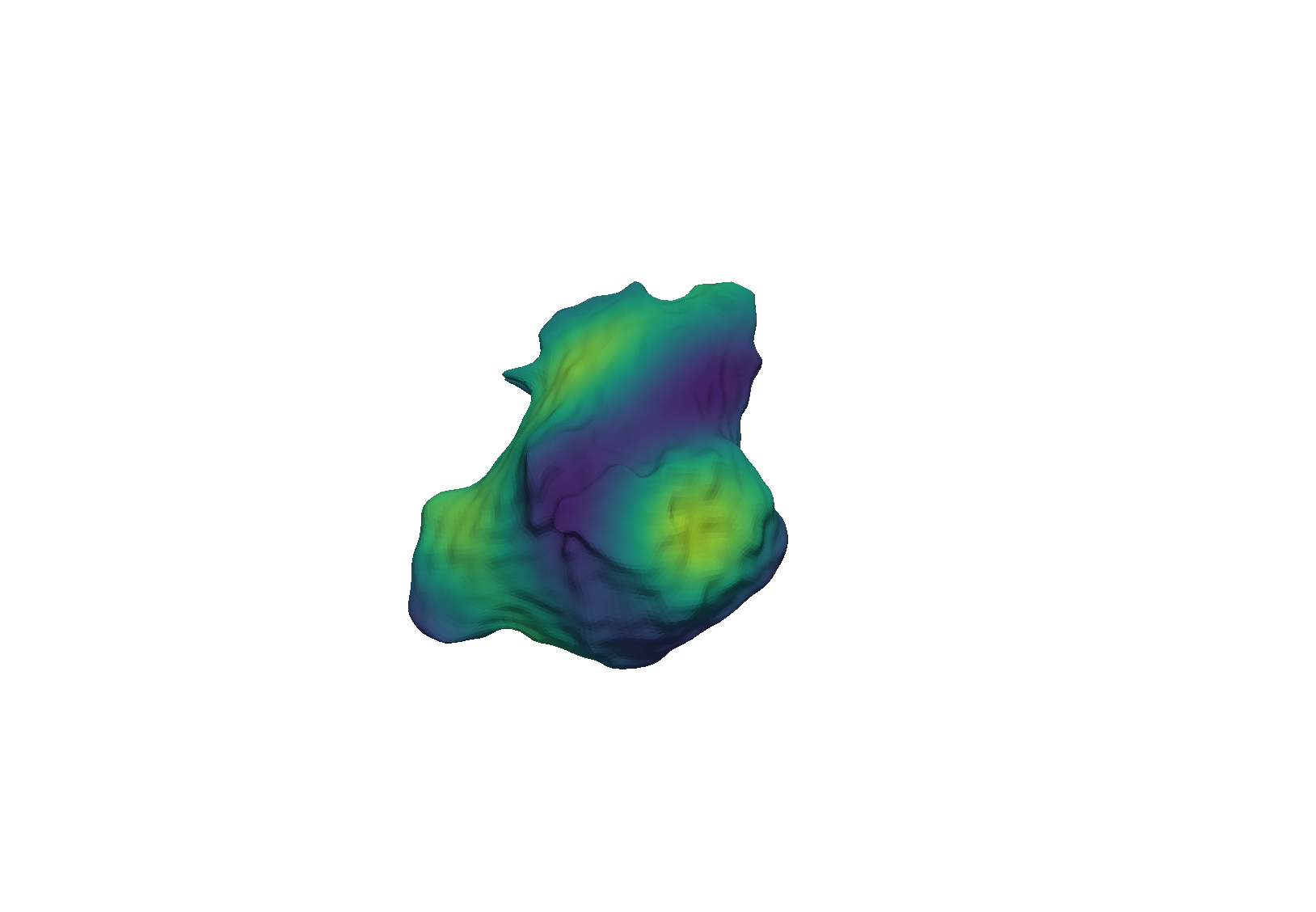}
\includegraphics[width=0.22\textwidth,clip=true,trim=15cm 10cm 15cm 10cm]{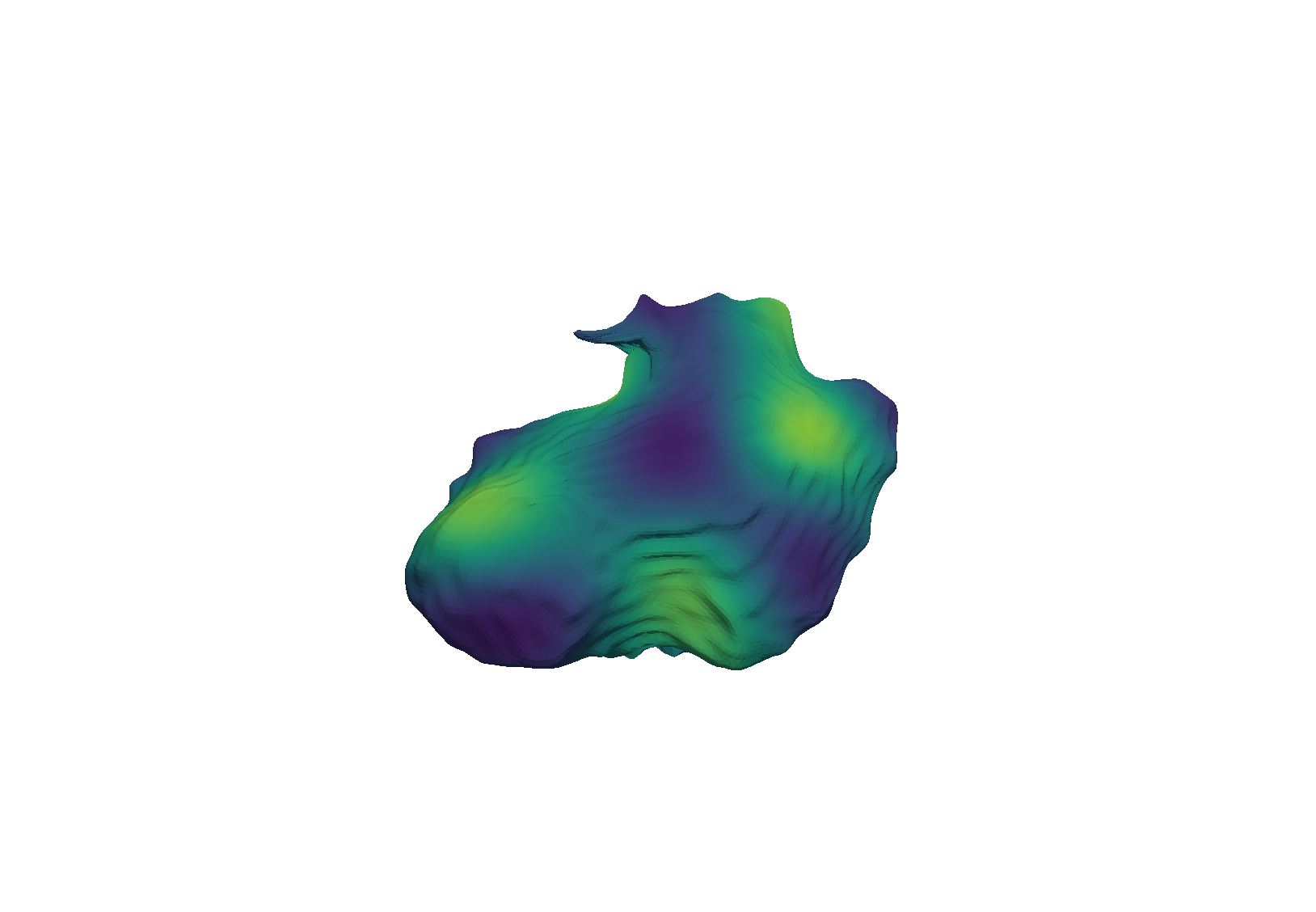}
\end{subfigure}

 \caption{Snapshots from results of Brusselator system on evolving surface at times $0,10,20,30,40,50,60,70$. See Section \ref{sec:pattern} for details on the simulation parameters.}
 \label{fig:bruss-evolving}
\end{figure}

The numerical scheme reads as before except that we have to solve for two fields now and the non-linear terms have been linearised using a first order Taylor expansion, i.e., for a function $f(u,w)$
\begin{multline*}
 f( u^{(m)}, w^{(m)} ) \approx f( u^{(m-1)}, w^{(m-1)} ) \\
 + f_u( u^{(m-1)}, w^{(m-1)} ) ( u^{(m)} - u^{(m-1)} ) + f_w( u^{(m-1)}, w^{(m-1)} ) ( w^{(m)} - w^{(m-1)} ).
\end{multline*}
To ensure stability, we take a computational time step $10^{-4}$ and interpolate between new meshes which are introduced every $1$ time unit until time $70$. We take a linear interpolation in time of the position of each node. Snapshots of the solution are shown in Figure~\ref{fig:bruss-evolving} and a video of all time steps is available at
\begin{quote}
  \url{http://www.personal.leeds.ac.uk/~scstr/vid/imagedata-video.html}
\end{quote}
We observe that the number and position of spots is changing over time. On stationary surfaces these are related to the eigenvalues of the Laplace-Beltrami operator which depend on the geometry. A careful assessment of the impact of the surface evolution is left for future studies.
%

%% file: source/conclusion_submitted.tex
\section{Conclusion}

A novel approach to solve reaction-advection-diffusion equations on moving surfaces which are obtained from time series of 3D cell image data was presented. 
By processing the image data a time series of triangulations of the same topology is obtained. An ALE finite element method then 
is set up to approximate solutions to the surface PDEs. 

The approach was used to simulate FRAP experiments with a recovery process based on diffusion along the membrane surface only. 
By fitting recovery curves a dependence on the  geometry of the photobleached region of interest was noticed. For experimentalists 
it is important to note that even in our idealised FRAP setup where we have assumed an idealised bleach region, which only extends 
through the dorsal side of the membrane and is perfectly symmetrical unlike real bleach profiles,  the computed half-times 
for recovery differ by factors up to 10.
For this purpose, the domain $\Gamma_b(t)$ to measure the recovery was obtained by simply intersecting the bleach profile (ROI) with the 
actual surface as this is the standard in the analysis of FRAP data \cite{ElzJan2009,GooKen2005}. This motivated restricting our numerical experiments  to small times in which the cell surface changes only a little.  However, over a longer time there is a significant evolution of the cell.
  This suggests that in a  
FRAP  model the bleached region  should be advected  with the material velocity in order to get more accurate recovery results.
This would motivate a change in the experimental set up so that the bleached region can be tracked.
Furthermore the analysis would benefit from a higher resolution data both in time and space which 
can be achieved by modern technology \cite{GaoShaCheBet14}. 

A constant diffusivity was assumed in our FRAP simulations. However, inhomogeneous diffusivities have been proposed, for instance, in 
dependence on the local geometry \cite{WeiBreAnd05} or on the composition of the cell membrane \cite{VasGhoGreEppFox04}. 
Inverse problems such as identifying the diffusivity of membrane resident proteins will require fluorescence data 
related to the protein density in addition to position data for the cell membrane. From such data, it may also be possible 
to infer information on the material velocity which is underlying the 
continuum approach to the surface PDEs. Of particular interest will be models for the lateral transport due to a gross 
motion of the membrane which often is modelled as a viscous fluid \cite{RAMB15}. 

In a second numerical experiment a pattern-forming reaction diffusion system on the cell surface was simulated. 
Also here an impact of the changing geometry, in this case in terms of the pattern, was observed. 

Our aim in this paper is  to point out that simulations of PDEs may be carried out on time evolving domains generated from data and 
that the results of standard models depend on the evolving geometry.  We expect that the capabilities of this computational surface finite 
element methodology may be used in parameter identification. For example by formulating  inverse problems to infer reaction rates or 
diffusivities of membrane resident species. But the results also suggest that
more systematic investigations of models and their dependence on geometry can  be carried out on realistic domains.  
For example, target models for investigation might be  reaction-diffusion systems for cell polarisation as are employed in models for 
cell motility \cite{NeiMacWeb11,EllStiVen12}. See also \cite{CroEllLad15} for a study on the quantification of such a model. 
Significant extensions might be required to  account for processes within and outside of the cell,\cite{MacMacNol16}, 
or raft formation in realistic geometries \cite{GarKamRaeRoe15}.